\documentclass[11pt]{article}
\usepackage{pxfonts}
\usepackage{yfonts}
\usepackage{dsfont}

\usepackage{graphicx}
\usepackage{relsize}

\def\bel{\begin{equation}\label}
\def\eeq{\end{equation}}
\def\ds{\displaystyle}

\def\mt{\longrightarrow}
\def\v{\vskip 1em}
\def\vsk{\vskip 40em}
\def\ve{\varepsilon}
\def\R{\mathds R}
\def\Z{\mathds Z}
\def\C{\mathfrak{B}}

\def\A{{\bf A}}
\def\B{{\bf B}}

\def\L{{\bf L}}
\def\U{{\bf U}}
\def\V{{\bf V}}

\def\I{{\bf I}}

\def\alpha{\alphaup}
\def\beta{\betaup}
\def\gamma{\gammaup}
\def\delta{\deltaup}
\def\xi{{\xiup}}
\def\eta{{\etaup}}
\def\tau{{\tauup}}
\def\rho{{\rhoup}}
\def\phi{{\phiup}}
\def\psi{{\psiup}}
\def\lambda{{\lambdaup}}
\def\omega{\omegaup}
\def\varphi{{\varphiup}}
\def\gamma{{\gammaup}}
\def\c{{\bf c}}

\def\({\left(}
\def\){\right)}

\newtheorem{remark}{Remark}[section]

\begin{document}
\[\begin{array}{cc}\ds\hbox{\LARGE{\bf Weighted $\L^p$-norm inequality of multi-parameter } }
\\\\\ds
\hbox{\LARGE{\bf fractional integration} }
\end{array}\]

\[\hbox{Chuhan Sun~~~~~and~~~~~Zipeng Wang}\]

\begin{abstract}~~~
We study a family of strong fractional integral operators whose kernels have  singularity  on every coordinate subspace. We prove a two-weight $\L^p$-norm inequality provided that the corresponding multi-parameter $\vartheta$-bump characteristic is finite.
\end{abstract}

\section{Introduction}
\setcounter{equation}{0}
Let $0<\alpha<n$ and define
\bel{I_alpha}
\I_\alpha f(x)~=~\int_{\R^n} f(y)\left({1\over |x-y|}\right)^{n-\alpha} dy.
\eeq
In 1928, Hardy and Littlewood \cite{Hardy-Littlewood} have first proved the $\L^p\mt\L^q$-boundedness of $\I_\alpha$
when $n=1$. Ten years later, Sobolev \cite{Sobolev} extended this result to the higher dimensional spaces. Today, it is well known as  Hardy-Littlewood-Sobolev inequality.   

Over the second half of 20th century, the regularity of $\I_\alpha$ has been extensively studied  for its   weighted norm analogue. A number of classical results have established,   for example by 
 Stein and Weiss \cite{Stein-Weiss}, Muckenhoupt and Wheeden \cite{Muckenhoupt-Wheeden}, 
Fefferman and Muckenhoupt \cite{Fefferman-Muckenhoupt}, Coifman and Fefferman \cite{Coifman-Fefferman}, Perez \cite{Perez} and Sawyer and Wheeden \cite{Sawyer-Wheeden}.

Consider $\omega^p$, $\sigma^{-{p\over p-1}}$ be non-negative, locally integrable. For $1<p\leq q<\infty$ and $\vartheta>1$, the $\vartheta$-bump characteristic is defined by 
\bel{A_pq alpha}
\begin{array}{lr}\ds
\A_{pq\vartheta}^{\alpha}(\omega,\sigma)
~=~\sup_{Q~\subset~\R^n}
|Q|^{{\alpha\over n}-{1\over p}+{1\over q}}\left\{{1\over |Q|}\int_{Q}\omega^{q\vartheta}(x)dx\right\}^{1\over q\vartheta} \left\{{1\over |Q|}\int_{Q}\left({1\over \sigma}\right)^{p\vartheta\over p-1}(x)dx\right\}^{p-1\over p\vartheta}
\end{array}
\eeq
where $Q$ is a cube in $\R^n$.  $\A_{pq\vartheta}^{\alpha}(\omega,\sigma)<\infty$ is also known as the Fefferman-Phong condition, initially introduced for  $p=q$. See \cite{Fefferman} for more background. 

$\diamond$ {\small Throughout, $\C>0$ is regarded as a generic constant depending on its subindices.}

{\bf Theorem A}~ ~( Sawyer and Wheeden, 1992 ) ~~{\it Let $\I_\alpha$ defined in (\ref{I_alpha}) for $0<\alpha<n$. We have
\bel{Sawyer-Wheeden Theorem}
\left\| \omega\I_{\alpha}f\right\|_{\L^q(\R^n)}~\leq~\C_{p~q~\vartheta~\alpha}~\A_{pq\vartheta}^{\alpha}(\omega,\sigma)~\left\| f\sigma\right\|_{\L^p(\R^n)},\qquad 1<p\leq q<\infty
\eeq
if $\A_{pq\vartheta}^{\alpha}(\omega,\sigma)<\infty$ for some $\vartheta>1$.}
\begin{remark}The constant $\C_{p~q~\vartheta~\alpha}~ \A_{pq\vartheta}^{\alpha}(\omega, \sigma)$ inside (\ref{Sawyer-Wheeden Theorem}) is not written explicitly in the original version of the statement. But, it can be computed directly by carrying out the proof given in section 2 of  Sawyer and Wheeden \cite{Sawyer-Wheeden}.
\end{remark}
We aim to give an extension of {\bf Theorem A} in product spaces by considering so-called the strong fractional integral operator. 
Let $0<\alpha<n$, $0<\beta<m$. Define
\bel{I_ab}
\I_{\alpha\beta}f(x,y)~=~\iint_{\R^n\times\R^m} f(u,v) \left({1\over |x-u|}\right)^{n-\alpha}\left({1\over |y-v|}\right)^{m-\beta}dudv.
\eeq
The study of certain operators  that  commute with a multi-parameter family of dilations,  dates back to the time of Jessen, Marcinkiewicz and Zygmund. Some important works were accomplished 
by Fefferman \cite{R.Fefferman}-\cite{R.Fefferman''}, Chang and Fefferman \cite{Chang-Fefferman},
Cordoba and Fefferman \cite{Cordoba-Fefferman},     Fefferman and Stein \cite{R.Fefferman-Stein},    M\"{u}ller, Ricci and Stein \cite{M.R.S},
Journ\'{e} \cite{Journe'} and Pipher \cite{Pipher}.

 For $1<p\leq q<\infty$ and $\vartheta>1$, we define
 \bel{A_pq}
\begin{array}{lr}\ds
\A_{pq\vartheta}^{\alpha\beta}(\omega,\sigma)
~=~\sup_{Q\times P~\subset~\R^n\times\R^m}
\\ \ds
|Q|^{{\alpha\over n}-{1\over p}+{1\over q}}|P|^{{\beta\over m}-{1\over p}+{1\over q}}\left\{{1\over |Q||P|}\iint_{Q\times P}\omega^{q\vartheta}(x,y)dxdy\right\}^{1\over q\vartheta} \left\{{1\over |Q||P|}\iint_{Q\times P}\left({1\over \sigma}\right)^{p\vartheta\over p-1}(x,y)dxdy\right\}^{p-1\over p\vartheta}
\end{array}
\eeq
where $Q$, $P$ are cubes in $\R^n$ and $\R^m$ respectively.
\v
{\bf Theorem A*}~~{\it Let $\I_{\alpha\beta}$ defined in (\ref{I_ab}) for $0<\alpha<n$, $0<\beta<m$. We have
\bel{Norm Ineq}
\left\| \omega\I_{\alpha\beta}f\right\|_{\L^q(\R^n\times\R^m)}~\leq~\C_{p~q~\vartheta~\alpha~\beta}~\A_{pq\vartheta}^{\alpha\beta}(\omega,\sigma)~\left\| f\sigma\right\|_{\L^p(\R^n\times\R^m)},\qquad 1<p\leq q<\infty
\eeq
if $\A_{pq\vartheta}^{\alpha\beta}(\omega,\sigma)<\infty$ for some $\vartheta>1$.}

{\bf Theorem A*} is proved for $1<p<q<\infty$ by Sawyer and Wang \cite{Sawyer-Wang}. The regarding method  requires $p$ strictly less than $q$. The critical endpoint case $p=q$ was left open.

\section{Statement of main result}
\setcounter{equation}{0}
{\bf Theorem One}~~{\it Let $\I_{\alpha\beta}$ defined in (\ref{I_ab}) for $0<\alpha<n$, $0<\beta<m$. We have
\bel{Norm Ineq}
\left\| \omega\I_{\alpha\beta}f\right\|_{\L^p(\R^n\times\R^m)}~\leq~\C_{p~\vartheta~\alpha~\beta}~\A_{p\vartheta}^{\alpha\beta}(\omega,\sigma)~\left\| f\sigma\right\|_{\L^p(\R^n\times\R^m)},\qquad 1<p<\infty
\eeq
if 
\bel{A sup}
\begin{array}{lr}\ds
\A_{p\vartheta}^{\alpha\beta}(\omega,\sigma)
~=~\sup_{Q\times P~\subset~\R^n\times\R^m}
\\ \ds
|Q|^{{\alpha\over n}-{1\over \vartheta}}|P|^{{\beta\over m}-{1\over \vartheta}}\left\{\iint_{Q\times P}\omega^{p\vartheta}(x,y)dxdy\right\}^{1\over p\vartheta} \left\{\iint_{Q\times P}\left({1\over \sigma}\right)^{p\vartheta\over p-1}(x,y)dxdy\right\}^{p-1\over p\vartheta}~<~\infty
\end{array}
\eeq
for some $\vartheta>1$.}

In order to prove {\bf Theorem One}, we investigate a new framework whereas $\R^n\times\R^m$ is decomposed into an infinitely many dyadic cones. Each partial operator defined on one of these dyadic cones is essentially an one-parameter fractional integral operator, satisfying the desired regularity. Furthermore, the operator's norm decays exponentially  as the cone's eccentricity  getting large.

\section{Cone decomposition on the product space}
\setcounter{equation}{0}
For every $\ell\in\Z$, we define
\bel{Partial}
\Delta_\ell \I_{\alpha\beta} f(x,y)~=~\iint_{\Gamma_\ell(x,y)} f(u,v) \left({1\over |x-u|}\right)^{n-\alpha}\left({1\over |y-v|}\right)^{m-\beta}dudv
\eeq
where
\bel{Cone}
\begin{array}{lr}\ds
\Gamma_\ell(x,y)~=~\left\{(u,v)\in\R^n\times\R^m \colon~2^{-\ell}\leq {|x-u|\over|y-v|}< 2^{-\ell+1}\right\}.
\end{array}
\eeq
\begin{figure}[h]
\centering
\includegraphics[scale=0.40]{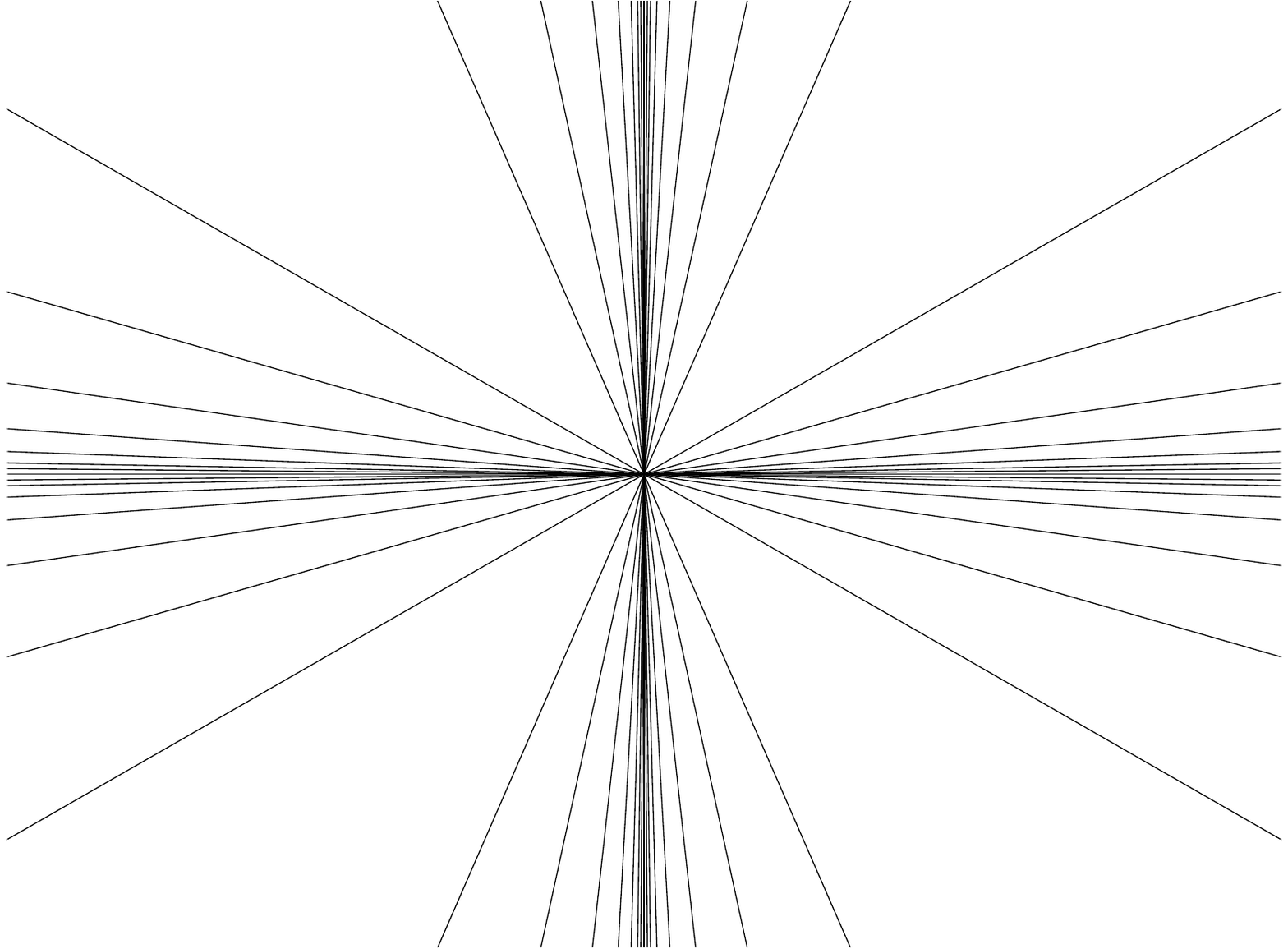}
\end{figure}

Recall (\ref{A sup}). By using H\"{o}lder inequality, we find $\A_{pt}^{\alpha\beta}(\omega,\sigma)\leq \A_{p\vartheta}^{\alpha\beta}(\omega,\sigma)$ for every $1<t\leq\vartheta$.

We define
\bel{A l}
\begin{array}{lr}\ds
\A_{pt}^{\alpha\beta}(\ell\colon\omega,\sigma)
~=~\sup_{Q\times P~\subset~\R^n\times\R^m~\colon~ |Q|^{1\over n}/|P|^{1\over m}=2^{-\ell}}
\\ \ds
|Q|^{{\alpha\over n}-{1\over t}}|P|^{{\beta\over m}-{1\over t}}\left\{\iint_{Q\times P}\omega^{pt}(x,y)dxdy\right\}^{1\over pt} \left\{\iint_{Q\times P}\left({1\over \sigma}\right)^{pt\over p-1}(x,y)dxdy\right\}^{p-1\over pt},\qquad \ell\in\Z.
\end{array}
\eeq
Given $\ell>0$, $2^{-\ell}Q$ is a dilate of $Q$ where
$|2^{-\ell}Q|^{1\over n}=2^{-\ell}|Q|^{1\over n}$. 
We have
\bel{A-Characteristic Dila Q}
\begin{array}{lr}\ds
|Q|^{{\alpha\over n}-{1\over t}}|P|^{{\beta\over m}-{1\over t}}\left\{\iint_{Q\times P}\omega^{pt}\left(2^{-\ell}x, y\right)dxdy\right\}^{1\over pt} \left\{\iint_{Q\times P}\left({1\over \sigma}\right)^{pt\over p-1}\left(2^{-\ell}x, y\right)dxdy\right\}^{p-1\over pt}
\\\\ \ds
~=~|Q|^{\alpha\over n}|P|^{\beta\over m}\left\{{1\over |Q||P|}\iint_{Q\times P}\omega^{pt}\left(2^{-\ell}x, y\right)dxdy\right\}^{1\over pt} \left\{{1\over |Q||P|}\iint_{Q\times P}\left({1\over \sigma}\right)^{pt\over p-1}\left(2^{-\ell}x, y\right)dxdy\right\}^{p-1\over pt}
\\\\ \ds
~=~2^{ \alpha\ell}|2^{-\ell}Q|^{\alpha\over n}|P|^{\beta\over m}\left\{{1\over|2^{-\ell}Q||P|}\iint_{2^{-\ell}Q\times P} \omega^{pt}\left( x,y\right)dxdy\right\}^{1\over pt}\left\{{1\over|2^{-\ell}Q||P|}\iint_{2^{-\ell}Q\times P} \left({1\over \sigma}\right)^{pt\over p-1}\left( x,y\right)dxdy\right\}^{p-1\over pt}
\\ \ds~~~~~~~~~~~~~~~~~~~~~~~~~~~~~~~~~~~~~~~~~~~~~~~~~~~~~~~~~~~~~~~~~~~~~~~~~~~~~~~~~~~~~~~~~~~~~~~~~~~~~~~~~~~~~~~~~~~~~~~~~~~~\hbox{\small{$x\mt 2^\ell x$}}
\\ \ds
~\leq~2^{ \alpha\ell}~\A_{pt}^{\alpha\beta}(\ell\colon\omega,\sigma)\qquad \hbox{\small{if~~ $|Q|^{1\over n}=|P|^{1\over m}$ ~~by (\ref{A l})}}.
\end{array}
\eeq
A similar estimate to (\ref{A-Characteristic Dila Q}) shows
\bel{A-Characteristic Dila P}
\begin{array}{lr}\ds
|Q|^{{\alpha\over n}-{1\over t}}|P|^{{\beta\over m}-{1\over t}}\left\{\iint_{P\times Q}\omega^{pt}\left( x,2^\ell y\right)dxdy\right\}^{1\over pt} \left\{\iint_{P\times Q}\left({1\over \sigma}\right)^{pt\over p-1}\left( x, 2^\ell y\right)dxdy\right\}^{p-1\over pt}
\\\\ \ds
~\leq~2^{ -\beta\ell}~\A_{pt}^{\alpha\beta}(\ell\colon\omega,\sigma),\qquad\ell\leq0,\qquad \hbox{\small{$|Q|^{1\over n}=|P|^{1\over m}$}}.
\end{array}
\eeq
Now, recall {\bf Theorem A} from section 1. We have
\bel{One-para Norm Ineq}
\begin{array}{lr}\ds
\left\{\iint_{\R^n\times\R^m}\left\{\iint_{\R^n\times\R^m} f(u,v) \left[{1\over \sqrt{|x-u|^2+|y-v|^2}}\right]^{n+m-(\alpha+\beta)}dudv\right\}^p\omega^p(x,y)dxdy\right\}^{1\over p}
\\\\ \ds
~\leq~\C_{p~t~\alpha~\beta}~ \A_{pt}^{\alpha+\beta}(\omega, \sigma) ~\left\{\iint_{\R^n\times\R^m} \Big(f\sigma\Big)^p(x,y)dxdy\right\}^{1\over p},\qquad 1<p<\infty
\end{array}
\eeq
 if 
\bel{A-Chara r-bump}
\begin{array}{lr}\ds
\A_{pt}^{\alpha+\beta}(\omega,\sigma)
~=~\sup_{Q\times P~\subset~\R^n\times\R^m~\colon~ |Q|^{1\over n}=|P|^{1\over m}}
\\ \ds
|Q|^{{\alpha\over n}-{1\over t}}|P|^{{\beta\over m}-{1\over t}}\left\{\iint_{Q\times P}\omega^{pt}(x,y)dxdy\right\}^{1\over pt} \left\{\iint_{Q\times P}\left({1\over \sigma}\right)^{pt\over p-1}(x,y)dxdy\right\}^{p-1\over pt}~<~\infty.
\end{array}
\eeq
By applying (\ref{One-para Norm Ineq})-(\ref{A-Chara r-bump})  and  using (\ref{A-Characteristic Dila Q})-(\ref{A-Characteristic Dila P}), we find
\bel{Regularity est Q}
\begin{array}{lr}\ds
\left\{\iint_{\R^n\times\R^m}\left\{\iint_{\R^n\times\R^m}f\left(2^{-\ell}u,v\right)\left[{1\over \sqrt{|x-u|^2+|y-v|^2}}\right]^{n+m-(\alpha+\beta)} dudv\right\}^p\omega^p\left(2^{-\ell} x,y\right)dxdy\right\}^{1\over p} 
\\\\ \ds
~\leq~\C_{p~t~\alpha~\beta}~2^{ \alpha\ell}~\A_{pt}^{\alpha\beta}(\ell\colon\omega,\sigma)\left\{\iint_{\R^n\times\R^m} \Big(f\sigma\Big)^p\left(2^{-\ell}x,y\right)dxdy\right\}^{1\over p},\qquad \ell>0
\end{array}
\eeq
and
\bel{Regularity est P}
\begin{array}{lr}\ds
\left\{\iint_{\R^n\times\R^m}\left\{\iint_{\R^n\times\R^m}f\left( u,2^\ell v\right)\left[{1\over \sqrt{|x-u|^2+|y-v|^2}}\right]^{n+m-(\alpha+\beta)} dudv\right\}^p\omega^p\left(  x,2^\ell y\right)dxdy\right\}^{1\over p} 
\\\\ \ds
~\leq~\C_{p~t~\alpha~\beta}~2^{-\beta\ell}~\A_{pt}^{\alpha\beta}(\ell\colon\omega,\sigma)\left\{\iint_{\R^n\times\R^m} \Big(f\sigma\Big)^p\left(x,2^\ell y\right)dxdy\right\}^{1\over p},\qquad \ell\leq0.
\end{array}
\eeq
Recall (\ref{Partial})-(\ref{Cone}). We have
\bel{Dilation Est Q}
\begin{array}{lr}\ds
\left\{\iint_{\R^n\times\R^m} \Big(\Delta_\ell\I_{\alpha\beta} f\Big)^p(x,y)\omega^p(x,y)dxdy\right\}^{1\over p}
\\\\ \ds
~=~\left\{\iint_{\R^n\times\R^m}\left\{ \iint_{\Gamma_\ell(x,y)}f(u,v)\left({1\over |x-u|}\right)^{n-\alpha}
\left({1\over |y-v|}\right)^{m-\beta}
du dv\right\}^p\omega^p(x,y)dxdy\right\}^{1\over p}
 \\\\ \ds
 ~=~\left\{\iint_{\R^n\times\R^m}\left\{\iint_{\Gamma_o(x,y)}f\left(2^{-\ell}u, v\right)\left({1\over 2^{-\ell}|x-u|}\right)^{n-\alpha} \left({1\over  |y-v|}\right)^{m-\beta}2^{-\ell n} dudv\right\}^p\omega^p\left(2^{-\ell}x, y\right)2^{-\ell n} dxdy\right\}^{1\over p}
 \\ \ds~~~~~~~~~~~~~~~~~~~~~~~~~~~~~~~~~~~~~~~~~~~~~~~~~~~~~~~~~~~~~
 ~~~~~~~~~~~~~~~~~~~~~~~~~~~~~~~~~~~~~~~~~~~~~~~~~~~~~~~~~~~~~~~~~
 \hbox{\small{$x\mt2^{-\ell}x$,~$u\mt 2^{-\ell}u$}}
 \\\\ \ds
 ~\leq~\C~2^{-\ell\big[\alpha+{n\over p}\big]}\left\{\iint_{\R^n\times\R^m}\left\{ \iint_{\R^n\times\R^m}f\left(2^{-\ell}u,  v\right)\left[{1\over \sqrt{|x-u|^2+|y-v|^2}}\right]^{n+m-(\alpha+\beta)}dudv\right\}^p\omega^p\left(2^{-\ell} x, y\right)dxdy\right\}^{1\over p}
 \\\\ \ds
 ~\leq~\C_{p~t~\alpha~\beta}~2^{-\ell\big[\alpha+{n\over p}\big]} 2^{\alpha\ell }\A_{pt}^{\alpha\beta}\left(\ell~\colon\omega,\sigma\right)\left\{ \iint_{\R^n\times\R^m} \Big(f\sigma\Big)^p\left(2^{-\ell} x, y\right)dxdy\right\}^{1\over p}\qquad \hbox{\small{by (\ref{Regularity est Q})}}
  \\\\ \ds
   ~=~\C_{p~t~\alpha~\beta}~\A_{pt}^{\alpha\beta}\left(\ell~\colon\omega,\sigma\right)2^{-\big({n\over p}\big)\ell}\left\{ \iint_{\R^n\times\R^m} \Big(f\sigma\Big)^p\left(x,y\right) 2^{\ell n}dxdy\right\}^{1\over p}\qquad \hbox{\small{$x\mt 2^\ell x$}}
 \\\\ \ds
 ~=~\C_{p~t~\alpha~\beta}~\A_{pt}^{\alpha\beta}\left(\ell~\colon\omega,\sigma\right)\left\{ \iint_{\R^n\times\R^m} \Big(f\sigma\Big)^p\left(x,y\right)dxdy\right\}^{1\over p},\qquad \ell>0.
\end{array}
\eeq
A similar estimate to (\ref{Dilation Est Q}) by using (\ref{Regularity est P}) instead of (\ref{Regularity est Q}) shows
\bel{Dilation Est P}
\begin{array}{lr}\ds
\left\{\iint_{\R^n\times\R^m} \Big(\Delta_\ell\I_{\alpha\beta} f\Big)^p(x,y)\omega^p(x,y)dxdy\right\}^{1\over p}
 \\\\ \ds
 ~\leq~\C_{p~t~\alpha~\beta}~\A_{pt}^{\alpha\beta}\left(\ell~\colon\omega,\sigma\right)\left\{ \iint_{\R^n\times\R^m} \Big(f\sigma\Big)^p\left(x,y\right)dxdy\right\}^{1\over p},\qquad \ell\leq0.
\end{array}
\eeq
From (\ref{Dilation Est Q})-(\ref{Dilation Est P}), we conclude
\bel{Regularity EST}
\begin{array}{lr}\ds
\left\| \omega\Delta_\ell\I_{\alpha\beta} f\right\|_{\L^p(\R^n\times\R^m)}
~\leq~\C_{p~t~\alpha~\beta}~\A_{pt}^{\alpha\beta}\left(\ell~\colon\omega,\sigma\right)~\left\| f\sigma\right\|_{\L^p(\R^n\times\R^m)},~~~~ 1<p<\infty,~~~~ \ell\in\Z.
\end{array}
\eeq
Note that 
\bel{Sum}
\I_{\alpha\beta}f(x,y)~=~\sum_{\ell\in\Z} \Delta_\ell \I_{\alpha\beta}f(x,y).
\eeq
We finish the proof of {\bf Theorem One} by applying Minkowski inequality together with the next result. 

{\bf Lemma One}~~{\it There exists an $\ve>0$ such that 
\bel{Regularity EST*}
\begin{array}{lr}\ds
\left\| \omega\Delta_\ell\I_{\alpha\beta} f\right\|_{\L^p(\R^n\times\R^m)}
~\leq~\C_{p~\vartheta~\alpha~\beta}~2^{-\ve|\ell|}~\A^{\alpha\beta}_{p\vartheta}(\omega,\sigma)~\left\| f\sigma\right\|_{\L^p(\R^n\times\R^m)},\qquad 1<p<\infty
\end{array}
\eeq
for every $\ell\in\Z$.}

\section{Proof of Lemma One}
\setcounter{equation}{0}
By symmetry, we prove (\ref{Regularity EST*}) for $\ell>0$ only. 
Let $1<r\leq\vartheta$ and denote
\bel{A, B}
\begin{array}{cc}\ds
\A_{pr}^{\alpha\beta}\Big[{Q\times P}\Big] (\omega,\sigma)~=~|Q|^{{\alpha\over n}-{1\over r}}|P|^{{\beta\over m}-{1\over r}} \B_{pr}\Big[{Q\times P}\Big] (\omega,\sigma),
\\\\ \ds
\B_{pr}\Big[{Q\times P}\Big] (\omega,\sigma)~=~ \left\{\iint_{Q\times P}\omega^{pr}(x,y)dxdy\right\}^{1\over pr} \left\{\iint_{Q\times P}\left({1\over \sigma}\right)^{pr\over p-1}(x,y)dxdy\right\}^{p-1\over pr}.
\end{array}
\eeq
Moreover,  $Q^\ell$ is a cube  having the same center of $Q$ with  $|Q^\ell|^{1\over n}=2^{-\ell}|Q|^{1\over n}$. 

Clearly, ${\B_{pr}\left[{Q^\ell\times P}\right] (\omega,\sigma)\over \B_{pr}[{Q\times P}] (\omega,\sigma)}\leq1$ for every $\ell>0$ and $Q\times P\subset\R^n\times\R^m$. 
We have
\bel{Compara Est Q <}
\begin{array}{lr}\ds
\A_{p\vartheta}^{\alpha\beta}\Big[{Q^\ell\times P}\Big] (\omega,\sigma)
~=~|Q^\ell|^{{\alpha\over n}-{1\over \vartheta}}|P|^{{\beta\over m}-{1\over \vartheta}}\B_{p\vartheta}\Big[{ Q^\ell\times P}\Big] (\omega,\sigma)
\\\\ \ds
~\leq~\C~|Q^\ell|^{{\alpha\over n}-{1\over \vartheta}}|P|^{{\beta\over m}-{1\over \vartheta}}\B_{p\vartheta}\Big[{ Q\times P}\Big] (\omega,\sigma)
\\\\ \ds
~=~\C~2^{-\ell\big[\alpha-{n\over \vartheta}\big]}|Q|^{{\alpha\over n}-{1\over \vartheta}}|P|^{{\beta\over m}-{1\over \vartheta}}\B_{p\vartheta}\Big[{Q\times P}\Big] (\omega,\sigma)
\\\\ \ds
~=~\C~2^{-\ell\big[\alpha-{n\over \vartheta}\big]} |Q|^{{\alpha\over n}-{1\over \vartheta}}|P|^{{\beta\over m}-{1\over \vartheta}} \left\{\iint_{Q\times P}\omega^{p\vartheta}(x,y)dxdy\right\}^{1\over p\vartheta} \left\{\iint_{Q\times P}\left({1\over \sigma}\right)^{p\vartheta\over p-1}(x,y)dxdy\right\}^{p-1\over p\vartheta}.
\end{array}
\eeq
Suppose $\alpha>{n\over \vartheta}$.
By taking the supremum over all $Q^\ell\times P \subset \R^n\times\R^m$ with $|Q|^{1\over n}=|P|^{1\over m}$ in (\ref{Compara Est Q <}), we conclude
\bel{EST >}
\begin{array}{lr}\ds
\A_{p\vartheta}^{\alpha\beta}\left(\ell\colon\omega,\sigma\right)~\leq~\C~2^{-\ell\big[\alpha-{n\over \vartheta}\big]} ~\A_{p\vartheta}^{\alpha\beta}\left(0\colon\omega,\sigma\right)
\\\\ \ds~~~~~~~~~~~~~~~~~~~~~
~\leq~\C~2^{-\ell\big[\alpha-{n\over \vartheta}\big]} ~\A_{p\vartheta}^{\alpha\beta}\left(\omega,\sigma\right).
\end{array}
\eeq
Note that (\ref{EST >}) together with (\ref{Regularity EST}) further imply (\ref{Regularity EST*}).

From now on, our assertion is on $\alpha\leq{n\over \vartheta}$. 

Given $Q\times P\subset\R^n\times\R^m$, 
we assume
\bel{Compara condition Hypothesis}
{\B_{pr}\left[{Q^\ell\times P}\right] (\omega,\sigma)\over \B_{pr}\Big[{Q\times P}\Big] (\omega,\sigma)}~\leq~\C~2^{\ell\big[\alpha-{n\over r}\big]},\qquad 1<r\leq\vartheta,\qquad\ell>0.
\eeq
Otherwise, there is a sequence  $\C_\ell, \ell=1,2,\ldots$ such that  $\C_\ell\mt\infty$ 
as $\ell\mt\infty$ and
\bel{Compara condition l sequence}
{\B_{pr}\left[{ Q^\ell\times P}\right] (\omega,\sigma)\over \B_{pr}\Big[{Q\times P}\Big] (\omega,\sigma)}~>~\C_\ell~2^{\ell\big[\alpha-{n\over r}\big]},\qquad r\in(1,\vartheta].
\eeq
This further implies
\bel{Compara Est Q}
\begin{array}{lr}\ds
\A_{pr}^{\alpha\beta}\Big[{Q^\ell\times P}\Big] (\omega,\sigma)
~=~|Q^\ell|^{{\alpha\over n}-{1\over r}}|P|^{{\beta\over m}-{1\over r}}\B_{pr}\Big[{ Q^\ell\times P}\Big] (\omega,\sigma)
\\\\ \ds
~=~2^{-\ell\big[\alpha-{n\over r}\big]}|Q|^{{\alpha\over n}-{1\over r}}|P|^{{\beta\over m}-{1\over r}}\B_{pr}\Big[{ Q^\ell\times P}\Big] (\omega,\sigma)
\\\\ \ds
~>~\C_\ell~2^{-\ell\big[\alpha-{n\over r}\big]}|Q|^{{\alpha\over n}-{1\over r}}|P|^{{\beta\over m}-{1\over r}}2^{\ell\big[\alpha-{n\over r}\big]}\B_{pr}\Big[{Q\times P}\Big] (\omega,\sigma)\qquad\hbox{\small{by (\ref{Compara condition l sequence})}}
\\\\ \ds
~=~\C_\ell~|Q|^{{\alpha\over n}-{1\over r}}|P|^{{\beta\over m}-{1\over r}} \left\{\iint_{Q\times P}\omega^{pr}(x,y)dxdy\right\}^{1\over pr} \left\{\iint_{Q\times P}\left({1\over \sigma}\right)^{pr\over p-1}(x,y)dxdy\right\}^{p-1\over pr}.
\end{array}
\eeq
Because $\A^{\alpha\beta}_{pr}(\omega, \sigma)\leq\A_{p\vartheta}^{\alpha\beta}(\omega,\sigma)<\infty$, it is clear that 
\bel{l bound}
\A_{pr}^{\alpha\beta}\Big[{Q^\ell\times P}\Big] (\omega,\sigma)~<~\infty\qquad \hbox{as}~~\ell\mt\infty.
\eeq
By putting together (\ref{Compara condition l sequence}), (\ref{Compara Est Q}) and (\ref{l bound}), we must have
\bel{A=0}
\begin{array}{rl}\ds
|Q|^{{\alpha\over n}-{1\over r}}|P|^{{\beta\over m}-{1\over r}}\left\{\iint_{Q\times P}\omega^{pr}(x,y)dxdy\right\}^{1\over pr} \left\{\iint_{Q\times P}\left({1\over \sigma}\right)^{pr\over p-1}(x,y)dxdy\right\}^{p-1\over pr}~=~0.
\end{array}
\eeq
This can only occurs if either $\omega=0$ or $\sigma^{-1}=0$ on $Q\times P\subset\R^n\times\R^m$. We shall exclude such trivial cases. 

Therefore, by using (\ref{A, B}), we rewrite (\ref{Compara condition Hypothesis}) as
\bel{Compara condition rewrite}
\begin{array}{lr}\ds
{\B_{pr}\left[{Q^\ell\times P}\right] (\omega,\sigma)\over \B_{pr}\Big[{Q\times P}\Big] (\omega,\sigma)}~=~
{\ds\left\{\iint_{Q^\ell\times P}\omega^{pr}(x,y)dxdy\right\}^{1\over pr}
\over \ds\left\{\iint_{Q\times P}\omega^{pr}(x,y)dxdy\right\}^{1\over pr}}
{ \ds\left\{\iint_{ Q^\ell\times P}\left({1\over \sigma}\right)^{pr\over p-1}(x,y)dxdy\right\}^{p-1\over pr}\over \ds \left\{\iint_{Q\times P}\left({1\over \sigma}\right)^{pr\over p-1}(x,y)dxdy\right\}^{p-1\over pr}}
\\\\ \ds~~~~~~~~~~~~~~~~~~~~~~~~~~~~~~~
~\leq~\C~2^{\ell\big[\alpha-{n\over r}\big]},\qquad 1<r\leq\vartheta,\qquad \ell>0.
\end{array}
\eeq
Furthermore,  consider
$1<r<s<\vartheta$  where $\alpha\leq{n\over \vartheta}<{n\over s}<{n\over r}$.

Because of (\ref{Compara condition rewrite}),   we can find $\lambda_\omega, \lambda_\sigma,\mu_\omega, \mu_\sigma\in[0,1]$
such that
\bel{RD} 
\begin{array}{cc}\ds
{\left\{\ds\iint_{ Q^\ell\times P}\omega^{pr}(x,y)dxdy\right\}^{1\over{pr}}\over\ds\left\{\iint_{Q\times P}\omega^{pr}(x,y)dxdy\right\}^{1\over{pr}}}~\leq ~\C~2^{\ell\big[\alpha-{n\over r}\big]\lambda_\omega},\qquad
{\ds\left\{\iint_{ Q^\ell\times P}\left(1\over\sigma\right)^{pr\over p-1}(x,y)dxdy\right\}^{p-1\over{pr}}\over\ds\left\{\iint_{Q\times P}\left(1\over\sigma\right)^{pr\over p-1}(x,y)dxdy\right\}^{p-1\over{pr}}}~\leq~ \C~ 2^{\ell\big[\alpha-{n\over r}\big]\lambda_\sigma},
\\\\ \ds
{\ds\left\{\iint_{Q^\ell\times P}\omega^{ps}(x,y)dxdy\right\}^{1\over{ps}}\over\ds\left\{\iint_{Q\times P}\omega^{ps}(x,y)dxdy\right\}^{1\over{ps}}}~\leq~ \C~ 2^{\ell\big[\alpha-{n\over s}\big]\mu_\omega},
\qquad
{\ds\left\{\iint_{Q^\ell\times P}\left(1\over\sigma\right)^{ps\over p-1}(x,y)dxdy\right\}^{p-1\over{ps}}\over\ds\left\{\iint_{Q\times P}\left(1\over\sigma\right)^{ps\over p-1}(x,y)dxdy\right\}^{p-1\over{ps}}}~\leq~ \C 2^{\ell\big[\alpha-{n\over s}\big]\mu_\sigma}\end{array}
\eeq
where
\bel{lambda mu functions}
\begin{array}{ccc}\ds
\lambda_\omega~=~\lambda_\omega\left( pr,  Q\times P,\ell\right),\qquad \lambda_\sigma~=~\lambda_\sigma\left( {pr\over p-1}, Q\times P,\ell\right),\qquad \lambda_\omega+\lambda_\sigma~\ge~1,
\\\\ \ds
\mu_\omega~=~\mu_\omega\left( ps, Q\times P,\ell\right),\qquad \mu_\sigma~=~\mu_\sigma\left( {ps\over p-1},Q\times P,\ell\right),\qquad \mu_\omega+\mu_\sigma~\ge~1.
\end{array}
\eeq
\begin{remark}
We keep in mind that $\lambda_\omega, \lambda_\sigma$ and $\mu_\omega, \mu_\sigma$  depend on $pr, ps, Q\times P, \ell$ and ${pr\over p-1}, {ps\over p-1},Q\times P, \ell$ respectively.
\end{remark}
Denote
\bel{gamma_omega}
\gamma_\omega~=~(\mu_\omega-\lambda_\omega)\left[\alpha-{n\over r}\right]+n\left({1\over r}-{1\over s}\right)\left(\mu_\omega-{1\over p}\right)
\eeq
and
\bel{gamma_sigma}
\gamma_\sigma~=~(\mu_\sigma-\lambda_\sigma)\left[\alpha-{n\over r}\right]+n\left({1\over r}-{1\over s}\right)\left(\mu_\sigma-{p-1\over p}\right).
\eeq
A direct computation shows
\bel{gamma omega computa}
\begin{array}{lr}\ds
(\lambda_\omega-\mu_\omega)\left[\alpha-{n\over s}\right]+n\left({1\over s}-{1\over r}\right)\left(\lambda_\omega-{1\over p}\right)
\\\\ \ds
~=~(\lambda_\omega-\mu_\omega)\left[\alpha-{n\over r}\right]+n\left({1\over s}-{1\over r}\right)\left(\mu_\omega-{1\over p}\right)~=~-\gamma_\omega
\end{array}
\eeq
and
\bel{gamma sigma computa}
\begin{array}{lr}\ds
(\lambda_\sigma-\mu_\sigma)\left[\alpha-{n\over s}\right]+n\left({1\over s}-{1\over r}\right)\left(\lambda_\sigma-{p-1\over p}\right)
\\\\ \ds
~=~(\lambda_\sigma-\mu_\sigma)\left[\alpha-{n\over r}\right]+n\left({1\over s}-{1\over r}\right)\left(\mu_\sigma-{p-1\over p}\right)~=~-\gamma_\sigma.
\end{array}
\eeq
Recall (\ref{A, B}). We have
\bel{Apr>omega}
\begin{array}{lr}\ds
\A_{pr}^{\alpha\beta}\Big[Q^\ell\times P\Big](\omega, \sigma)
\\\\ \ds
~=~|Q^\ell|^{\alpha\over n}|P|^{\beta\over m}\left\{{1\over|Q^\ell||P|}\iint_{Q^\ell\times P}\omega^{pr}(x,y)dxdy\right\}^{1\over{pr}}\left\{{1\over |Q^\ell||P|}\iint_{Q^\ell\times P}{\left(1\over\sigma\right)}^{pr\over p-1}(x,y)dxdy\right\}^{p-1\over pr}
\\\\ \ds
~\leq~|Q^\ell|^{\alpha\over n}|P|^{\beta\over m}\left\{{1\over |Q^\ell||P|}\iint_{Q^\ell\times P}\omega^{pr}(x,y)dxdy\right\}^{1\over{pr}}\left\{{1\over |Q^\ell||P|}\iint_{Q^\ell\times P}{\left(1\over\sigma\right)}^{ps\over p-1}(x,y)dxdy\right\}^{p-1\over ps}
\\ \ds
~~~~~~~~~~~~~~~~~~~~~~~~~~~~~~~~~~~~~~~~~~~~~~~~~~~~~~~~~~~~~~~~~~~~~~~~~~~~~~~~~~~~~~~~~~~~~~~~~\hbox{ \small{by H\"{o}lder  inequality}}
\\\\ \ds
~=~|Q^\ell|^{{\alpha\over n}-{1\over pr}-{p-1\over ps}}|P|^{{\beta\over m}-{1\over pr}-{p-1\over ps}}\left\{\iint_{Q^\ell\times P}\omega^{pr}(x,y)dxdy\right\}^{1\over{pr}}\left\{\iint_{Q^\ell\times P}{\left(1\over\sigma\right)}^{ps\over p-1}(x,y)dxdy\right\}^{p-1\over ps}
\\\\ \ds
~\leq~\C~ |Q|^{{\alpha\over n}-{1\over pr}-{p-1\over ps}}|P|^{{\beta\over m}-{1\over pr}-{p-1\over ps}}2^{-\ell\Big[\alpha-\frac{n}{pr}-\frac{n(p-1)}{ps}\Big]}2^{\ell\big[\alpha-{n\over r}\big]\lambda_\omega}2^{\ell\big[\alpha-{n\over s}\big]\mu_\sigma}
\\ \ds~~~~~~~
\left\{\iint_{Q\times P}\omega^{pr}(x,y)dxdy\right\}^{1\over{pr}}\left\{\iint_{Q\times P}{\left(1\over\sigma\right)}^{ps\over p-1}(x,y)dxdy\right\}^{p-1\over ps}\qquad\hbox{\small{by (\ref{RD})}}
\\\\ \ds
~\leq~\C~ |Q|^{{\alpha\over n}-{1\over pr}-{p-1\over ps}}|P|^{{\beta\over m}-{1\over pr}-{p-1\over ps}}2^{-\ell\Big[\alpha-\frac{n}{pr}-\frac{n(p-1)}{ps}\Big]}2^{\ell\big[\alpha-{n\over r}\big]\lambda_\omega}2^{\ell\big[\alpha-{n\over s}\big](1-\mu_\omega)}
\\ \ds~~~~~~~
\left\{\iint_{Q\times P}\omega^{pr}(x,y)dxdy\right\}^{1\over{pr}}\left\{\iint_{Q\times P}{\left(1\over\sigma\right)}^{ps\over p-1}(x,y)dxdy\right\}^{p-1\over ps}\qquad \hbox{\small{( $\mu_\omega+\mu_\sigma\ge1$ )}}
\\\\ \ds
~=~\C~2^{-\ell\Big[(\mu_\omega-\lambda_\omega)\big[\alpha-{n\over r}\big]+n\big({1\over r}-{1\over s}\big)\big(\mu_\omega-{1\over p}\big)\Big]}
\\ \ds
~~~~~~~
|Q|^{\alpha\over n}|P|^{\beta\over m}\left\{{1\over|Q||P|}\iint_{Q\times P}\omega^{pr}(x,y)dxdy\right\}^{1\over{pr}}\left\{{1\over|Q| |P|}\iint_{Q\times P}{\left(1\over\sigma\right)}^{ps\over p-1}(x,y)dxdy\right\}^{p-1\over ps}
\\\\ \ds
~\leq~\C~2^{-\ell\Big[(\mu_\omega-\lambda_\omega)\big[\alpha-{n\over r}\big]+n\big({1\over r}-{1\over s}\big)\big(\mu_\omega-{1\over p}\big)\Big]}
\\ \ds
~~~~~~~|Q|^{\alpha\over n}|P|^{\beta\over n}\left\{{1\over|Q| |P|}\iint_{Q\times P}\omega^{ps}(x,y)dxdy\right\}^{1\over{ps}}\left\{{1\over|Q| |P|}\iint_{Q\times P}{\left(1\over\sigma\right)}^{ps\over p-1}(x,y)dxdy\right\}^{p-1\over ps}
\\ \ds
~~~~~~~~~~~~~~~~~~~~~~~~~~~~~~~~~~~~~~~~~~~~~~~~~~~~~~~~~~~~~~~~~~~~~~~~~~~~~~~~~~~~~~~~~~ \hbox{ \small{by H\"{o}lder  inequality}}
\\\\ \ds
~=~\C~2^{-\ell\gamma_\omega}
\A_{ps}^{\alpha\beta}\Big[Q\times P\Big](\omega, \sigma)
\qquad
\hbox{\small{by (\ref{gamma_omega})}}.
\end{array}
\eeq
By carrying out a similar estimate to (\ref{Apr>omega}) with $\lambda_\omega$ replaced by $1-\lambda_\sigma$ and using (\ref{gamma_sigma}) instead of (\ref{gamma_omega}), we obtain
\bel{EST gamma_sigma >}
\begin{array}{rl}\ds
\A_{pr}^{\alpha\beta}\left[Q^\ell\times P\right](\omega,\sigma)
~\leq~\C~2^{-\ell\gamma_\sigma}~\A_{ps}^{\alpha\beta}\left[Q\times P\right](\omega,\sigma).
\end{array}
\eeq

On the other hand, we  have
\bel{Apr<omega}
\begin{array}{lr}\ds
\A_{pr}^{\alpha\beta}\Big[Q^\ell\times P\Big](\omega, \sigma)
\\\\ \ds
~=~|Q^\ell|^{\alpha\over n}|P|^{\beta\over m}\left\{{1\over |Q^\ell||P|}\iint_{Q^\ell\times P}\omega^{pr}(x,y)dxdy\right\}^{1\over{pr}}\left\{{1\over|Q^\ell||P|}\iint_{Q^\ell\times P}{\left(1\over\sigma\right)}^{pr\over p-1}(x,y)dxdy\right\}^{p-1\over pr}
\\\\ \ds
~\leq~|Q^\ell|^{\alpha\over n}|P|^{\beta\over m}\left\{{1\over|Q^\ell||P|}\iint_{Q^\ell\times P}\omega^{ps}(x,y)dxdy\right\}^{1\over{ps}}\left\{{1\over |Q^\ell||P|}\iint_{Q^\ell\times P}{\left(1\over\sigma\right)}^{pr\over p-1}(x,y)dxdy\right\}^{p-1\over pr}
\\ \ds
~~~~~~~~~~~~~~~~~~~~~~~~~~~~~~~~~~~~~~~~~~~~~~~~~~~~~~~~~~~~~~~~~~~~~~~~~~~~~~~~~~~~~~~~~~~~~~~ \hbox{ \small{by H\"{o}lder  inequality}}
\\\\ \ds
~=~|Q^\ell|^{{\alpha\over n}-{1\over ps}-{p-1\over pr}}|P|^{{\beta\over m}-{1\over ps}-{p-1\over pr}}\left\{\iint_{Q^\ell\times P}\omega^{ps}(x,y)dxdy\right\}^{1\over{ps}}\left\{\iint_{Q^\ell\times P}{\left(1\over\sigma\right)}^{pr\over p-1}(x,y)dxdy\right\}^{p-1\over pr}
\\\\ \ds
~\leq~\C~ |Q|^{{\alpha\over n}-{1\over ps}-{p-1\over pr}}|P|^{{\beta\over m}-{1\over ps}-{p-1\over pr}}2^{-\ell\Big[\alpha-\frac{n}{ps}-\frac{n(p-1)}{pr}\Big]} 2^{\ell\big[\alpha-{n\over s}\big]\mu_\omega}2^{\ell\big[\alpha-{n\over r}\big]\lambda_\sigma}
\\ \ds
~~~~~~~\left\{\iint_{Q\times P}\omega^{ps}(x,y)dxdy\right\}^{1\over{ps}}\left\{\iint_{Q\times P}{\left(1\over\sigma\right)}^{pr\over p-1}(x,y)dxdy\right\}^{p-1\over pr}\qquad\hbox{\small{by (\ref{RD})}}
\\\\ \ds
~\leq~\C~ |Q|^{{\alpha\over n}-{1\over ps}-{p-1\over pr}}|P|^{{\beta\over m}-{1\over ps}-{p-1\over pr}}2^{-\ell\Big[\alpha-\frac{n}{ps}-\frac{n(p-1)}{pr}\Big]} 2^{\ell\big[\alpha-{n\over s}\big]\mu_\omega}2^{\ell\big[\alpha-{n\over r}\big](1-\lambda_\omega)}
\\ \ds
~~~~~~~\left\{\iint_{Q\times P}\omega^{ps}(x,y)dxdy\right\}^{1\over{ps}}\left\{\iint_{Q\times P}{\left(1\over\sigma\right)}^{pr\over p-1}(x,y)dxdy\right\}^{p-1\over pr} \qquad \hbox{\small{( $\lambda_\omega+\lambda_\sigma\ge1$ )}}

\\\\ \ds
~=~\C~2^{-\ell\Big[(\lambda_\omega-\mu_\omega)\big[\alpha-{n\over s}\big]+n\big({1\over s}-{1\over r}\big)\big(\lambda_\omega-{1\over p}\big)\Big]}
\\ \ds
~~~~~~~
|Q|^{\alpha\over n}|P|^{\beta\over m}\left\{{1\over|Q| |P|}\iint_{Q\times P}\omega^{ps}(x,y)dxdy\right\}^{1\over{ps}}\left\{{1\over|Q| |P|}\iint_{Q\times P}{\left(1\over\sigma\right)}^{pr\over p-1}(x,y)dxdy\right\}^{p-1\over pr}
\\\\ \ds
~\leq~\C~2^{-\ell\Big[(\lambda_\omega-\mu_\omega)\big[\alpha-{n\over s}\big]+n\big({1\over s}-{1\over r}\big)\big(\lambda_\omega-{1\over p}\big)\Big]}
\\ \ds
~~~~~~~
|Q|^{\alpha\over n}|P|^{\beta\over m}\left\{{1\over|Q| |P|}\iint_{Q\times P}\omega^{ps}(x,y)dxdy\right\}^{1\over{ps}}\left\{{1\over|Q| |P|}\iint_{Q\times P}{\left(1\over\sigma\right)}^{ps\over p-1}(x,y)dxdy\right\}^{p-1\over ps}
\\ \ds
~~~~~~~~~~~~~~~~~~~~~~~~~~~~~~~~~~~~~~~~~~~~~~~~~~~~~~~~~~~~~~~~~~~~~~~~~~~~~~~~~~~~~~~~~~~ \hbox{ \small{by H\"{o}lder  inequality}}

\\\\ \ds
~=~\C~2^{\ell\gamma_\omega}
\A_{ps}^{\alpha\beta}\Big[Q\times P\Big](\omega, \sigma)
\qquad\hbox{\small{by (\ref{gamma omega computa})}}.
\end{array}
\eeq
By carrying out a similar estimate to (\ref{Apr<omega}) with $\lambda_\omega$ replaced by $1-\lambda_\sigma$ and using (\ref{gamma sigma computa}) instead of (\ref{gamma omega computa}), we obtain
\bel{EST gamma_sigma <}
\begin{array}{rl}\ds
\A_{pr}^{\alpha\beta}\left[Q^\ell\times P\right](\omega,\sigma)
~\leq~\C~2^{\ell\gamma_\sigma}~\A_{ps}^{\alpha\beta}\left[Q\times P\right](\omega,\sigma).
\end{array}
\eeq
\vsk

Recall (\ref{lambda mu functions}) and (\ref{gamma_omega})-(\ref{gamma_sigma}). We have $\lambda_\omega=\lambda_\omega(pr, Q\times P, \ell)$, $\mu_\omega= \mu_\omega(ps, Q\times P, \ell)$ and $\lambda_\sigma= \lambda_\sigma\left({pr\over p-1},Q\times P,\ell\right)$, $\mu_\sigma=\mu_\sigma\left({ps\over p-1},Q\times P, \ell\right)$. Hence that $\gamma_\omega=\gamma_\omega(pr,ps, Q\times P,\ell)$ and $\gamma_\sigma=\gamma_\sigma\left({pr\over p-1},{ps\over p-1}, Q\times P, \ell\right)$.

Let $1<r=r_*<\vartheta$  fixed.
Suppose   $\exists ~s\colon r_*<s<\vartheta$ such that $|\gamma_\omega|>\ve$ or $|\gamma_\sigma|>\ve$ for some $\ve>0$. From  (\ref{Apr>omega}) to (\ref{EST gamma_sigma <}), we conclude
\bel{EST gamma non-zero}
\begin{array}{lr}\ds
\A_{pr_*}^{\alpha\beta}\left[Q^\ell\times P\right](\omega,\sigma)~\leq~\C~2^{-\ell\max\big\{|\gamma_\omega|,~|\gamma_\sigma|\big\}}~\A_{ps}^{\alpha\beta}\left[Q\times P\right](\omega,\sigma)
\\\\ \ds~~~~~~~~~~~~~~~~~~~~~~~~~~~~~~~~
~\leq~\C~2^{-\ell\ve}~\A_{ps}^{\alpha\beta}(0\colon\omega,\sigma)\qquad \hbox{\small{if~ $|Q|^{1\over n}=|P|^{1\over m}$}}
\\\\ \ds~~~~~~~~~~~~~~~~~~~~~~~~~~~~~~~~
~\leq~\C~2^{-\ell\ve}~\A_{p\vartheta}^{\alpha\beta}(\omega,\sigma).
\end{array}
\eeq
By solving the equations in (\ref{gamma_omega}) and (\ref{gamma_sigma}), we find
\bel{mu function}
\begin{array}{lr}\ds
\mu_\omega(ps, Q\times P, \ell)~=~ {\ds \lambda_\omega(pr_*, Q\times P,\ell)\left[\alpha-{n\over r_*}\right]+{1\over p}\left({n\over r_*}-{n\over s}\right)\over \ds \alpha-{n\over s}}+{\gamma_\omega(pr_*,ps, Q\times P,\ell)\over\ds \alpha-{n\over s}},
\\ \ds~~~~~~~~~~~~~~~~~~~~~~~~~~~~~~~~~~~~~~~~~~~~~~~~~~~~~~~~~~~
~~~~~~~~~~~~~~~~~~~~~~~~~~~~~~~~~~~~~~~~~~~~~~~
~~~~~~~~~~~r_*<s<\vartheta
\end{array}
\eeq
and
\bel{sigma function}
\begin{array}{lr}\ds
\mu_\sigma\left({ps\over p-1},Q\times P, \ell\right)~=~ {\ds \lambda_\sigma\left({pr_*\over p-1},Q\times P,\ell\right)\left[\alpha-{n\over r_*}\right]+{p-1\over p}\left({n\over r_*}-{n\over s}\right)\over \ds \alpha-{n\over s}}
\\\\ \ds~~~~~~~~~~~~~~~~~~~~~~~~~~~~~~~~~~
+{\ds\gamma_\sigma\left({pr_*\over p-1},{ps\over p-1}, Q\times P, \ell\right)\over\ds \alpha-{n\over s}},
\qquad r_*<s<\vartheta.
\end{array}
\eeq
\begin{remark} We have
$\A_{p_1\eta}^{\alpha\beta}(\omega,\sigma)\leq\A_{p\vartheta}^{\alpha\beta}(\omega,\sigma)<\infty$  and $\A_{p_2\eta}^{\alpha\beta}(\omega,\sigma)\leq\A_{p\vartheta}^{\alpha\beta}(\omega,\sigma)$ for $1<\eta<\vartheta$ where 
${p\vartheta\over p-1}={p_1 \eta\over p_1-1}$ and $p\vartheta=p_2 \eta$. \end{remark}
Indeed, for every $Q\times P\subset\R^n\times\R^m$, 
\bel{p_i Est}
\begin{array}{lr}\ds
\A_{p\vartheta}^{\alpha\beta}\left[ Q\times P\right](\omega, \sigma)
\\\\ \ds
~=~|Q|^{\alpha\over n}|P|^{\beta\over m}\left\{{1\over |Q||P|}\iint_{Q\times P}\omega^{p\vartheta}(x,y)dxdy\right\}^{1\over p\vartheta} \left\{{1\over |Q||P|}\iint_{Q\times P}\left({1\over \sigma}\right)^{p_1\eta\over p_1-1}(x,y)dxdy\right\}^{p_1-1\over p_1\eta}
\\\\ \ds
~=~|Q|^{\alpha\over n}|P|^{\beta\over m}\left\{{1\over |Q||P|}\iint_{Q\times P}\omega^{p_2\eta}(x,y)dxdy\right\}^{1\over p_2\eta} \left\{{1\over |Q||P|}\iint_{Q\times P}\left({1\over \sigma}\right)^{p\vartheta\over p-1}(x,y)dxdy\right\}^{p-1\over p\vartheta}.
\end{array}
\eeq
{\bf Remark 4.2} is a consequence of using H\"{o}lder inequality in (\ref{p_i Est}) whereas $p_1<p<p_2$.

Therefore, we can repeat all regarding estimates from (\ref{A, B}) to (\ref{EST gamma_sigma <}) for every $1<r<s<\eta$ $w.r.t$ $p_i$ for $i=1,2$.
In compare to (\ref{gamma_omega})-(\ref{gamma_sigma}), we find
\bel{gamma_i, omega}
\begin{array}{lr}\ds
\gamma_{\omega~i}~=~\gamma_{\omega}\left(p_ir, p_i s, Q\times P,\ell\right)~=~
\\\\ \ds
\Big(\mu_\omega(p_is, Q\times P,\ell)-\lambda_\omega(p_ir, Q\times P, \ell)\Big)\left[\alpha-{n\over r}\right]+n\left({1\over r}-{1\over s}\right)\left(\mu_\omega(p_is,Q\times P,\ell)-{1\over p_i}\right),
~~
 i=1,2
\end{array}
\eeq
and
\bel{gamma_i, sigma}
\begin{array}{lr}\ds
\gamma_{\sigma~i}~=~\gamma_{\sigma}\left({p_ir\over p_i-1}, {p_i s\over p_i-1}, Q\times P,\ell\right)~=~
\\\\ \ds
\left[\mu_\sigma\left({p_is\over p_i-1},Q\times P,\ell\right)-\lambda_\sigma\left({p_ir\over p_i-1}, Q\times P,\ell\right)\right]\left[\alpha-{n\over r}\right]
+n\left({1\over r}-{1\over s}\right)\left(\mu_\sigma\left({p_is\over p_i-1},Q\times P,\ell\right)-{p_i-1\over p_i}\right),
\\ \ds~~~~~~~~~~~~~~~~~~~~~~~~~~~~~~~~~~~~~~~~~~~~~~~~~~~~~~~~~~~~~~~~~~
~~~~~~~~~~~~~~~~~~~~~~~~~~~~~~~~~~~~~~~~~~~~~~~~~~~~~~~~~~~~~~
i=1,2.
\end{array}
\eeq
Recall (\ref{EST gamma non-zero}). Let $1<r=r_*<\vartheta$  fixed. For $\eta$  close enough to $\vartheta$ in {\bf Remark 4.2}, we have $r_*<\eta$.
Suppose $\exists~ s\colon r_*<s<\eta$ such that $|\gamma_{\omega~i}|>\ve$ or $|\gamma_{\sigma~i}|>\ve$ for some $\ve>0$. We obtain
 \bel{EST gamma non-zero i}
\begin{array}{lr}\ds
\A_{p_ir_*}^{\alpha\beta}\left[Q^\ell\times P\right](\omega,\sigma)~\leq~\C~2^{-\ell\max\big\{|\gamma_{\omega~i}|,~|\gamma_{\sigma~i}|\big\}}~\A_{p_is}^{\alpha\beta}\left[Q\times P\right](\omega, \sigma)
\\\\ \ds~~~~~~~~~~~~~~~~~~~~~~~~~~~~~~~~
~\leq~\C~2^{-\ell\ve}~\A_{p_is}^{\alpha\beta}(0\colon\omega,\sigma)\qquad\hbox{\small{if~ $|Q|^{1\over n}=|P|^{1\over m}$}}
\\\\ \ds~~~~~~~~~~~~~~~~~~~~~~~~~~~~~~~~
~\leq~\C~2^{-\ell\ve}~\A_{p_i\eta}^{\alpha\beta}(\omega,\sigma)
\\\\ \ds~~~~~~~~~~~~~~~~~~~~~~~~~~~~~~~~
~\leq~\C~2^{-\ell\ve}~\A_{p\vartheta}^{\alpha\beta}(\omega,\sigma),\qquad i=1,2\qquad\hbox{\small{by {\bf Remark 4.2}}}.
\end{array}
\eeq
By solving the equation in (\ref{gamma_i, omega}) with $r=r_*$, we find
\bel{mu i function}
\begin{array}{lr}\ds
\mu_\omega(p_is, Q\times P,\ell)~=~ {\ds \lambda_\omega(p_ir_*, Q\times P,\ell)\left[\alpha-{n\over r_*}\right]+{1\over p_i}\left({n\over r_*}-{n\over s}\right)\over \ds \alpha-{n\over s}}+{\gamma_\omega(p_ir_*,p_is, Q\times P, \ell)\over\ds \alpha-{n\over s}},
\\ \ds~~~~~~~~~~~~~~~~~~~~~~~~~~~~~~~~~~~~~~~~~~~~~~~~~~~~~~~~~~~~~~~~~~~~~
~~~~~~~~~~~~~~~~~~~~~~~~~~
 r_*<s<\eta,\qquad i=1,2.
\end{array}
\eeq
Choose $\c_i>0$ so that $p_i=p\c_i$ for $i=1,2$. (\ref{mu i function}) can be written as
\bel{mu i function c}
\begin{array}{lr}\ds
\mu_\omega(p\c_i s,Q\times P,\ell)~=~ {\ds \lambda_\omega(p \c_i r_*, Q\times P,\ell)\left[\alpha-{n\over r_*}\right]+{1\over p\c_i}\left({n\over r_*}-{n\over s}\right)\over \ds \alpha-{n\over s}}+{\gamma_\omega(p\c_i r_*,p\c_i s, Q\times P, \ell)\over\ds \alpha-{n\over s}},
\\ \ds~~~~~~~~~~~~~~~~~~~~~~~~~~~~~~~~~~~~~~~~~~~~~~~~~~~~~~~~~~~~~~~~~~
~~~~~~~~~~~~~~~~~~~~~~~~~~~
 r_*<s<\eta,\qquad i=1,2.
\end{array}
\eeq
On the other hand, recall (\ref{mu function}). By changing variables $s\mt \c_i s$, we also have
\bel{mu function c}
\begin{array}{lr}\ds
\mu_\omega(p\c_i s, Q\times P,\ell)~=~ {\ds \lambda_\omega(pr_*, Q\times P, \ell)\left[\alpha-{n\over r_*}\right]+{1\over p}\left({n\over r_*}-{n\over \c_i s}\right)\over \ds \alpha-{n\over \c_i s}}+{\gamma_\omega(p r_*,p\c_i s, Q\times P,\ell)\over\ds \alpha-{n\over \c_i s}},
\\ \ds~~~~~~~~~~~~~~~~~~~~~~~~~~~~~~~~~~~~~~~~~~~~~~~~~~~~~~~~~~~~~~~~~~~~
~~~~~~~~~~~~~~~~~~~~~~
 r_*\c^{-1}_i<s< \vartheta \c^{-1}_i,\qquad i=1,2.
\end{array}
\eeq
By subtracting  (\ref{mu function c}) from (\ref{mu i function c}), we find
\bel{mu function c range}
\begin{array}{lr}\ds
{\ds \lambda_\omega(p \c_i r_*, Q\times P, \ell)\left[\alpha-{n\over r_*}\right]+{1\over p\c_i}\left({n\over r_*}-{n\over s}\right)\over \ds \alpha-{n\over s}}
~-~ {\ds \lambda_\omega(pr_*, Q\times P, \ell)\left[\alpha-{n\over r_*}\right]+{1\over p}\left({n\over r_*}-{n\over \c_i s}\right)\over \ds \alpha-{n\over \c_i s}}
\\\\ \ds
~=~{\gamma_\omega(p r_*,p\c_i s, Q\times P, \ell)\over\ds \alpha-{n\over \c_i s}}-{\gamma_\omega(p\c_i r_*,p\c_i s, Q\times P, \ell)\over\ds \alpha-{n\over s}},
\qquad
\max\left\{r_*, r_*\c^{-1}_i\right\}<s< \min\left\{\eta, \vartheta\c^{-1}_i\right\},
\\ \ds~~~~~~~~~~~~~~~~~~~~~~~~~~~~~~~~~~~~~~~~~~~~~~~~~~~~~~~~~~~~~~~~~~~~~~~~~~~~~~~~~~~~~~~~~~~~~~~~~~~~~~~~~~~~~~~~~~~
 i=1,2.
\end{array}
\eeq
For $\eta$ sufficiently close to $\vartheta$ in {\bf Remark 4.2}, we have  $\c_i,i=1,2$  close enough to $1$ such that $\max\left\{r_*, r_*\c^{-1}_i\right\}< \min\left\{\eta, \vartheta\c^{-1}_i\right\},~i=1,2$.

By multiplying $\left[\alpha-{n\over  s}\right]\left[\alpha-{n\over \c_i s}\right]$ to both sides of the equation in  (\ref{mu function c range}), we have
\bel{mu function c range rewrite}
\begin{array}{lr}\ds
\left[\alpha-{n\over \c_i s}\right] \left[\lambda_\omega(p \c_i r_*, Q\times P, \ell)\left[\alpha-{n\over r_*}\right]+{1\over p\c_i}\left({n\over r_*}-{n\over s}\right)\right]
\\\\ \ds
~-~ \left[\alpha-{n\over s}\right] \left[\lambda_\omega(pr_*, Q\times P, \ell)\left[\alpha-{n\over r_*}\right]+{1\over p}\left({n\over r_*}-{n\over \c_i s}\right)\right]
\\\\ \ds
~=~\left[\alpha-{n\over  s}\right]
\gamma_\omega(p r_*,p\c_i s, Q\times P, \ell)-\left[\alpha-{n\over \c_i s}\right]\gamma_\omega(p\c_i r_*,p\c_i s, Q\times P, \ell).
\end{array}
\eeq
The left hand side of (\ref{mu function c range rewrite}) equals
 \bel{Polynomial s^-1}
 {n^2\over p}\left[{1\over \c_i^2}-{1\over \c_i}\right] \left({1\over s}\right)^2+\U(\omega, p, \c_i, r_*,\alpha, Q\times P, \ell)\left({1\over s}\right)+\V(\omega, p, \c_i, r_*,\alpha, Q\times P, \ell)
\eeq
where
\bel{UV}
\begin{array}{lr}\ds
 \U(\omega, p, \c_i, r_*,\alpha, Q\times P, \ell)~=~
 \\\\ \ds
 \left[{\alpha\over p\c_i}+\lambda_\omega(pr_*, Q\times P, \ell)\left(\alpha-{n\over r_*}\right)+{n\over pr_*}\right]n
 - \left[{\alpha\over p}+\lambda_\omega(p\c_i r_*, Q\times P, \ell)\left(\alpha-{n\over r_*}\right)+{n\over p\c_ir_*}\right]{n\over \c_i}, 
\\\\ \ds
 \V(\omega, p, \c_i, r_*,\alpha,Q\times P, \ell)~=~
 \\\\ \ds
 \alpha\left[\lambda_\omega(p\c_i r_*, Q\times P, \ell)\left(\alpha-{n\over r_*}\right)+{n\over p\c_ir_*}\right]
 - \alpha\left[\lambda_\omega(pr_*, Q\times P, \ell)\left(\alpha-{n\over r_*}\right)+{n\over pr_*}\right].
 \end{array}
 \eeq
 Observe that, 
 for $\c_i\neq1, i=1,2$, (\ref{Polynomial s^-1}) 
  is a second degree polynomial of $s^{-1}$, regardless of $\lambda_\omega(p \c_i r_*, Q\times P, \ell)$ and $\lambda_\omega(pr_*, Q\times P, \ell)$ in (\ref{UV}). Therefore, the right hand side of (\ref{mu function c range rewrite}) does not  vanish uniformly.
 \begin{remark}
  We cannot find a convergent subsequence of $Q_j\times P_j\subset\R^n\times\R^m, \ell_j>0$, $j=1,2\ldots$ such that  $\gamma_\omega(p r_*,p\c_i s, Q_j\times P_j, \ell_j)\mt0$ and $\gamma_\omega(p\c_i r_*,p\c_i s, Q_j\times P_j, \ell_j)\mt0$  as $j\mt\infty$ for every 
    $\max\left\{r_*, r_*\c^{-1}_i\right\}<s< \min\left\{\eta, \vartheta\c^{-1}_i\right\}$. 
  \end{remark}
Consequently, there is an $\ve>0$ such that for every $Q\times P\subset\R^n\times\R^m$ and $\ell>0$, we either have 
\bel{non-zero omega}
\begin{array}{cc}\ds
\left|\gamma_\omega(pr_*, p\c_i s, Q\times P, \ell)\right|~>~\ve\qquad \hbox{for some ~~$r_*\c_i^{-1}<s<\vartheta\c_i^{-1}$}
\\ \ds
\Longleftrightarrow
\\ \ds
\left|\gamma_\omega(pr_*, ps, Q\times P, \ell)\right|~>~\ve\qquad\hbox{for some~~ $r_*<s<\vartheta$}
\end{array}
\eeq
 or
\bel{non-zero i omega}
\left|\gamma_\omega(p_ir_*,p_is, Q\times P, \ell)\right|~>~\ve,\qquad i=1,2,\qquad\hbox{for some~~ $r_*<s<\eta$.}
\eeq
A similar argument to (\ref{mu i function})-(\ref{UV}) with $\lambda_\sigma, \mu_\sigma$  instead of $\lambda_\omega, \mu_\omega$ shows that we either have 
\bel{non-zero sigma}
\left|\gamma_\sigma\left({pr_*\over p-1}, {ps\over p-1}, Q\times P, \ell\right)\right|~>~\ve\qquad\hbox{for some ~~$r_*<s<\vartheta$}
\eeq
or 
\bel{non-zero i sigma}
\left| \gamma_{\sigma}\left({p_ir_*\over p_i-1}, {p_i s\over p_i-1}, Q\times P, \ell\right)\right|~>~\ve,\qquad i=1,2,\qquad\hbox{for some ~~$r_*<s< \eta$.}
\eeq
Note that (\ref{non-zero omega}) and (\ref{non-zero sigma}) imply (\ref{EST gamma non-zero}).  (\ref{non-zero i omega}) and (\ref{non-zero i sigma}) imply (\ref{EST gamma non-zero i}). 
By using (\ref{EST gamma non-zero i}),  we have
\bel{A product}
\begin{array}{lr}\ds
\C~2^{-2\ell\ve}~\Big[\A_{p\vartheta}^{\alpha\beta}(\omega,\sigma)\Big]^2~\ge~
\prod_{i=1,2}\A_{p_i r_*}^{\alpha\beta}\left[Q^\ell\times P\right](\omega,\sigma)
\\\\ \ds
~\ge~\left\{|Q^\ell|^{{\alpha\over n}}|P|^{{\beta\over m}}\left\{{1\over|Q^\ell||P|}\iint_{Q^\ell\times P}\omega^{p_1r_*}(x,y)dxdy\right\}^{1\over p_1r_*} \left\{{1\over|Q^\ell||P|}\iint_{Q^\ell\times P}\left({1\over \sigma}\right)^{p_2r_*\over p_2-1}(x,y)dxdy\right\}^{p_2-1\over p_2r_*}\right\}^2
\\ \ds~~~~~~~~~~~~~~~~~~~~~~~~~~~~~~~~~~~~~~~~~~~~~~~~~~~~~~~~~~~~~~~~~~~~
~~~~~~~~~~~~~~~
\hbox{\small{by H\"{o}lder inequality and $p_1<p<p_2$}}.
\end{array}
\eeq
\begin{remark} As $\eta$ close enough to $\vartheta$, 
 we have $p_1<p<p_2$ sufficiently close for which there is a $1<t<r_*<s<\eta<\vartheta$ such that
$pt\leq p_1 r_*$ and ${pt\over p-1}\leq{p_2 r_*\over p_2-1}$.
\end{remark}
From (\ref{A product}) and {\bf Remark 4.4}, we have
\bel{EST eta}
\begin{array}{lr}\ds
\A_{pt}^{\alpha\beta}\left[Q^\ell\times P\right](\omega,\sigma)
\\\\ \ds
~\leq~|Q^\ell|^{{\alpha\over n}}|P|^{{\beta\over m}}\left\{{1\over|Q^\ell||P|}\iint_{Q^\ell\times P}\omega^{p_1r_*}(x,y)dxdy\right\}^{1\over p_1r_*} \left\{{1\over|Q^\ell||P|}\iint_{Q^\ell\times P}\left({1\over \sigma}\right)^{p_2r_*\over p_2-1}(x,y)dxdy\right\}^{p_2-1\over p_2r_*}
\\ \ds~~~~~~~~~~~~~~~~~~~~~~~~~~~~~~~~~~~~~~~~~~~~~~~~~~~~~~~~~~~~~
~~~~~~~~~~~~~~~~~~~~~~~~~~~~~~~~~~~~~~~~~~~~
\hbox{\small{by H\"{o}lder inequality}}
\\ \ds
~\leq~\C~2^{-\ve\ell}~\A_{p\vartheta}^{\alpha\beta}(\omega,\sigma).
\end{array}
\eeq
Moreover, $\A_{pt}^{\alpha\beta}\left[Q^\ell\times P\right](\omega,\sigma)\leq \A_{pr_*}^{\alpha\beta}\left[Q^\ell\times P\right](\omega,\sigma)$ as $t<r_*$.

Lastly,  (\ref{EST gamma non-zero}) and (\ref{EST eta}) hold for every $Q\times P\subset\R^n\times\R^m$ 
 with $|Q|^{1\over n}=|P|^{1\over m}$ and every $\ell>0$.  We conclude
  \bel{EST decay}
\begin{array}{lr}\ds
\A_{pt}^{\alpha\beta}(\ell\colon\omega, \sigma)
~\leq~\C~2^{-\ve\ell}~\A_{p\vartheta}^{\alpha\beta}(\omega, \sigma),\qquad \ell>0
\end{array}
\eeq
for some $\ve>0$ and $1<t<\vartheta$.
This together with (\ref{Regularity EST}) further imply (\ref{Regularity EST*}).

\end{document}